\documentclass[a4paper,11pt]{article}

\usepackage{latexsym}
\usepackage{amsfonts}
\usepackage{amsmath}
\usepackage{graphicx}

\newcommand{\takis}{\usefont{OMS}{ztmcm}{m}{n}}
\newcommand{\nicecal}[1]{\mbox{\takis #1}}

\newcommand{\Prob}{{\mathbb P}}

\newcommand{\Exp}{{\mathbb E}}

\newcommand{\nat}{{\mathbb N}}

\newcommand{\dd}{{\mathrm{d}}}
\newcommand{\e}{{\mathrm{e}}}

\newcommand{\oh}{{\mathrm{o}}}
\newcommand{\Oh}{{\mathrm{O}}}

\newcommand{\refs}[1]{{\rm (\ref{#1})}}
\newcommand{\convdistr}
{\stackrel{{\scriptsize\cal D}}{\rightarrow}}

\newcommand{\eqdistr}{\stackrel{{\scriptsize\cal D}}{=}}
\newcommand{\process}[1]{\left\{ {#1} \right\}}
\newcommand{\proces}[1]{\bigl\{ {#1} \bigr\}}
\newcommand{\proof}{{\it Proof}.\ }
\newcommand{\halmos}{{\mbox{\, \vspace{3mm}}} \hfill
\mbox{$\Box$}}
\newtheorem{Th}{Theorem}[section]
\newtheorem{Lemma}{Lemma}[section]

\newtheorem{Prop}{Proposition}[section]
\newtheorem{Cor}{Corollary}[section]

\newcommand{\logsim}{\approx_{\rm log}}
\renewcommand{\tau}{X}

\newcommand{\Fb}{\overline F}
\newcommand{\Gb}{\overline G}
\newcommand{\Hb}{\overline H}
\newcommand{\restart}{RESTART}

\begin{document}

\begin{center}
{\Large
{\bf Asymptotic Behavior of Total Times For Jobs\\
That Must Start Over If a Failure Occurs.}}\\[.2in]
S\o ren Asmussen\footnote{Department of Mathematical Sciences,
Aarhus University, Ny Munkegade, DK-8000 Aarhus C,
Denmark}, Pierre Fiorini\footnote{%
Department of Computer Science,
University of Southern Maine,
Portland, Maine, USA},
Lester Lipsky\footnote{
Department of Computer Science and Engineering,
University of Connecticut, Storrs, CT 06269-2155, USA}\\
Tomasz Rolski\footnote{Mathematical Institute, Wroclaw University, 
50-384 Wroclaw, Poland\\
Partially upported by a Marie Curie Transfer of Knowledge Fellowship:
Programme HAHAP MTKO-CT-2004-13389}, \&
Robert Sheahan${}^3$\\[.2in]
\today\\[.3in]
\end{center}
\begin{abstract}
Many processes must complete in the presence of failures.
Different systems respond to task failure in different ways. The system may
resume a failed task from the failure point (or a saved checkpoint shortly
before the failure point), it may give up on
the task and select a replacement task from the ready queue, or
it may restart the task.  The behavior of systems under the first two 
scenarios is well documented, 
but the third ({\em RESTART}) has resisted detailed
analysis.  In this paper we derive tight asymptotic relations between the 
distribution of
{\em task times} without failures to the {\em total time} when including 
failures, for any failure distribution.  In particular, we show that 
if the  task time distribution has an unbounded support
then the total time distribution $H$ is always heavy-tailed.
Asymptotic expressions are given for the tail of $H$
in various scenarios. The key ingredients of the analysis are the
Cram\'er--Lundberg asymptotics for geometric sums
and integral asymptotics, that in some cases are obtained
via Tauberian theorems and in some cases by bare-hand calculations.
\\[2mm]
\end{abstract}
\newpage
\section{Introduction}
For many systems failure is rare enough that it can be ignored, or dealt
with as an afterthought. For other systems, failure is common enough that the
design choice of how to deal with it may have a significant impact on the
performance of the system.  Consider a job that ordinarily would take
a time $T$ to be executed on some system (e.g., CPU).  
If at some time $U<T$ the processor 
fails, the job may take a {\em total time} $X\ge T$ to complete.
We let $F,G$ be the distributions of $T,U$ and $H=H_{F,G}$ the distribution
of $X$ which in addition to $F,G$ depends on the failure recovery scheme.

Many papers discuss methods of failure recovery
and analyze their complexity in one or more metrics, like {\em restartable
processors} in Chlebus {\em et al.}~\cite{Chlebus01}, 
or {\em stage checkpointing} in De~Prisco {\em et al.}~\cite{DePrisco94}, etc.
There are many specific and distinct failure recovery schemes, but they can
be grouped into three broad classes:
\begin{quote}
{\em RESUME}, also referred to as preemptive resume (prs);\\
{\em REPLACE}, also referred to as preemptive repeat different (prd);\\
{\em RESTART}, also referred to as preemptive repeat identical (pri).
\end{quote}
The analysis of the distribution function $H(x)=\Prob(X\le x)$ 
when the policy is {\em RESUME} or {\em REPLACE} was carried out by 
Kulkarni {\em et al.}~\cite{Kulkarni86}, \cite{Kul}.  
In the {\em RESUME} scenario, if there is
a processor failure while a job is being executed, after repair is 
implemented  the job can continue where it left off.  All that is required
mathematically is to remember the state of the system when failure occurred.
If repair time is an issue then the number of failures before final completion
must also be considered.  
In what follows, we ignore the time for repairs, with the
knowledge that this can be properly handled separately.
In the {\em REPLACE} situation, if a job fails, it is replaced by a different
job from the same distribution.  Here, no details concerning the previous job
are necessary in order to continue. 

The work by Kulkarni {\em et al.}~\cite{Kulkarni86}, \cite{Kul}, 
and Bobbio \& Trivedi~\cite{Bobbio90} clearly suggests that
if $F$ is phase-type or, more generally, matrix-exponential
(\cite{Lipsky}, \cite{RP}, \cite{APQ}), 
and $\Gb(u)=\Prob(U>u)=\e^{-\beta u}$,
then $H$ for the {\em RESUME} and {\em REPLACE} 
policies can also be represented by  matrix-exponential distributions. 
This means that they
could be analyzed entirely within a Markov chain framework.  

However, the \restart\ policy has resisted detailed analysis.  
The total time distribution $H$ under this policy 
was defined and examined through its Laplace transform in Kulkarni {\em et al.}
\cite{Kulkarni86}, \cite{Kul}.  They were able to show
that it definitely was {\em not} matrix-exponential, 
i.e., the Laplace transform cannot be rational, 
and therefore it cannot be solved in the
Markov Chain framework.  However, by numerically taking the
inverse Laplace transform (see Jagerman~\cite{Jager82}), 
Chimento \& Trivedi \cite{Chim93}
(following a model proposed by Castillo~\cite{Castillo80}) 
were able to find the \restart\
time distribution for a few cases, for a limited range of the total time 
($x\le 3\Exp T$).  The method seems to be unstable for larger $x$.  
It is this problem that interests us here.

There are many examples of where the \restart\ scenario is relevant.
The obvious one alluded to above involves execution of a program on
some computer.  If the computer fails, and the intermediate results are
not saved externally (e.g., by {\em checkpointing}), then the job must 
restart from the beginning.  As another example, one might wish to 
copy a file from a remote system using 
some standard protocol as FTP or HTTP.
The time it takes to copy a file is proportional
to its length.  A transmission error
immediately aborts the copy and discards the partially received data,
forcing the user to restart the copy from the beginning. 
Yet another example would be receiving `customer service'  by telephone.
Often, while dealing with a particular service agent, the connection
is broken.  Then the customer must redial the service center, and invariably
(after waiting in a queue) end up talking to a different agent, and have
to explain everything from the beginning.  

In our previous paper (Sheahan {\em et al.}~\cite{Sheahan06}), 
we derived an expression for the Laplace transform 
of the total time distribution $H$ for the \restart\ policy with
exponential failure rate, $\beta$.  We used 
it to get an expression for the moments $\Exp X^\ell$ of the total time.  
>From this we were able to argue that if the task-time distribution has an 
exponential tail, then $X$ has infinite moments for 
$\ell\ge\alpha=\lambda/\beta$, where $\lambda$ is the rate of the 
exponential tail (i.e., $\Fb(t)\sim c\e^{-\lambda t}$).  
This in turn implies that roughly $\Hb(x)\approx c/x^\alpha$,
i.e., $X$ is {\em power-tailed}.

This can have important implications, particularly in applications
where the time to finish a task is bounded by necessity.  
If a task takes too
long to complete it must be aborted, and an alternate solution provided.
In such applications it may be important to know $\Hb(x)$, for that
is the probability that a job will be aborted.
Power tails and {\em heavy tails}, generally, have a small but non-negligible 
probability of lasting for many, many times 
the mean, and thus $\Hb(x)$ for large $x$ can be important.  

In this paper we derive the
asymptotic behavior of $\Hb(x)$ as $x\to\infty$ under more general assumptions
than in \cite{Sheahan06} and in sharper form in a number of
important cases. As a first guess, one could believe that the
heaviness of $\Hb$ is determined by the heaviness of $\Fb$ and/or
$\Gb$. However, it turns out that the important feature is rather
how close are $\Fb$ and $\Gb$. This is demonstrated in a striking way
by the following result for the diagonal case:
\begin{Prop}\label{Prop02.01a} If $F=G$, then
$\displaystyle\overline H(x)\,\sim\,\frac{1}{\mu x}$.
\end{Prop}
Here $\mu=1/\Exp U$; we assume throughout in the paper that
$\mu>0$ and, for convenience, that $F,G$ have densities $f,g$
(this assumption can be relaxed at many places but we will not give
the details).
It is notable that no other conditions are required for
Proposition~\ref{Prop02.01a}, in particular no
precise information on how heavy the
common tail $\overline F=\overline G$ is!

The assumption that the task time distribution $F$ and the failure
time distribution $G$ be identical of course lacks interpretation
in the \restart\ setting. Thus, Proposition~\ref{Prop02.01a}
is more of a curiosity, which is further illustrated by the fact that a proof
can be given which is far simpler than the our proofs
for more general situations (see Section~\ref{S:Concl}). 
Nevertheless, the result indicates that the
tail behaviour of $H$ depends on  a delicate balance between the tails of
$F$ and $G$. We will also see that making $\overline F$ heavier 
makes  $\overline H$ heavier, making $\overline G$ heavier 
makes  $\overline H$ lighter. However, except for the
case when $F$ has a finite support, $\overline H$ is always heavy-tailed:

\begin{Prop}\label{Prop15.6a} Assume that the support of $F$
is unbounded.
Then\\
$\e^{\epsilon x}\overline H(x)\,\to\,\infty$ for any $\epsilon>0$.
\end{Prop}

In general, we will be able to obtain sharp asymptotics for $\Hb(x)$ when $\Fb$
and $\Gb$ are not too far away. The form of the result
(Theorem~\ref{Prop12.3a}) is regular variation of $\Hb$. For example, the
following result covers Gamma distributions:

\begin{Cor}\label{Cor23.5a} Assume $f,g$ belong to the class of densities
of asymptotic form $c\e^{-\lambda t}t^{\alpha-1}$,
with parameters $\lambda_F,\alpha_F,c_F$ for $f$ and  
$\lambda_G,\alpha_G,c_G$ for $g$.
Then $\overline H(x)\ \sim\ c_H\log^{\alpha_F-\alpha_H\alpha_G}x/
x^{\alpha_H}$, where $\alpha_H=
\lambda_F/\lambda_G$ and 
$$c_H\ =\ \frac{c_F\Gamma(\alpha_H)\lambda_G^{\alpha_H-1-\alpha_F
+\alpha_H\alpha_G}}{\mu^{\alpha_H}c_G^{\alpha_H}}
\,.$$
\end{Cor}
Numerical illustrations are given in \cite{Sheahan06} for $\alpha_G=1$
(i.e., $G$ exponential) and show
an excellent fit.

When $\Fb$ and $\Gb$ are more different (say $F$ has a power tail
and $G$ is exponential), we will derive logarithmic asymptotics for
$\Hb$. We will see forms
varying from extremely heavy tails like $1/\log^\alpha x$ over power tails
$1/x^\alpha$ to moderately
heavy tails like the Weibull tail $\e^{-x^\beta}$ with $\beta<1$.

The proofs of the paper are based on the representation
\begin{equation}\label{16.5a}X\ =\ T\,+\,S \ \ 
\mbox{where}\ \ N\ =\ \inf\process{n:\, U_{n+1}>T}\,,\ \ S\ =\  
\sum_{i=1}^{N}U_i\,,
\end{equation}
and $U_1,U_2,\ldots$ are the succesive failure times
(assumed i.i.d.\ with distribution $G$ and independent of $T$).
More precisely, we will use that given $T=t$, $S(t)=\sum_1^{N}U_i$ is a
compound geometric sum for which exponential Cram\'er-Lundberg
tail asymptotics is available, and uncondition to get our final results.
In Section~\ref{S:Main} we state
our main results, except for the case of a bounded task time $T$
which is treated in Section~\ref{S:Geom}. The analysis there departs
from a careful study of the case $T\equiv t$. 
Section~\ref{S:First} is devoted to the proof of the
following lemma, which is the key to the unbounded case:
\begin{Lemma}\label{L20.10c} Let $\mu=1/\Exp U$ and define
$$I_\pm(x,\epsilon)\ =\ \int_0^\infty 
\exp\process{-\mu\overline G(t)x(1\pm\epsilon)}
f(t)\,\dd t$$
Then for each $\epsilon>0$,
$$1-\epsilon\ \le\ \liminf_{x\to\infty}\frac{\overline H(x)}{I_+(x,\epsilon)}
\ \le\  \limsup_{x\to\infty}
\frac{\overline H(x)}{I_-(x,\epsilon)}\ \le\ 1+\epsilon \,.$$
\end{Lemma} 
This lemma essentially reduces
the  investigation of the asymptotics of $\overline H(x)$ to
the (not always straightforward!) purely analytical study of the asymptotics
of $I_+(x,\epsilon)$ and $I_-(x,\epsilon)$. Indeed, we will see
in Section~\ref{S:Further}
that once this is done, one is most often able 
to obtain the logarithmic asymptotics
of $\overline H(x)$ by letting $\epsilon\downarrow 0$, and in some cases
even the sharp asymptotics. Finally, Section~\ref{S:Concl} contains
some concluding remarks.

\section{Statement of Main Results} \label{S:Main}
Except for Proposition~\ref{Prop12.10b}, we will assume throughout
the paper that the support of $F$ is infinite.

We shall use the concept of {\em logarithmic
asymptotics} familiar from large deviations theory and
write $f(t)\,\logsim\,g(t)$ for two functions $f,g>0$ with
limits 0 at $t=\infty$ if $\log f(t)/\log g(t)\to 1$ as $t\to\infty$. We then
consider the following distribution classes:

$$\nicecal{F}_1:\ f(t)\,\logsim\,\e^{-\alpha t^\eta},\ \ \
\nicecal{F}_2:\ f(t)\,\logsim\,\frac{1}{t^{\alpha+1}},$$
$$\nicecal{G}_1:\ \overline G(t)\,\logsim\,\e^{-\beta t^\gamma},\ \ \
\nicecal{G}_2:\ \overline G(t)\,\logsim\,\frac{1}{t^\beta}
$$ 
Note that these definitions do not completely identify the tail behaviour of
$F,G$. For example, if $\overline G(t)\,\sim\,ct^\alpha \e^{-\beta t^\gamma}$,
then $G\in\nicecal{G}_1$, but one cannot identify $c,\alpha$, and if
$f(t)\,\sim\, c\log^\beta t/t^{\alpha+1}$, then $F\in\nicecal{F}_1$, 
but one cannot identify $c,\beta$.

Note also that  $f\in\nicecal{F}_2$ implies that 
$\overline F(t)\,\logsim\,{1}/{t^{\alpha}}$, and that a sufficient
(but not necessary) condition for $G\in\nicecal{G}_2$ is that
$g(t)$ is regularly varying with index $-\beta-1$.

Similarly to the definition of $f\logsim g$, we will write $f\approx_{\log\log}
g$ if  
$$\frac{\log \bigl(-\log f(t)\bigr)}{\log \bigl(-\log g(t)\bigr)}
\ \to \ 1\,.$$
See further  part (1:2) of Theorem \ref{Th22.10a} and 
Remark~E) in Section~\ref{S:Concl}.

With these distribution classes, we obtain a complete description of
the logarithmic asymptotics of $\overline H(x)$ 
except for the case $f\in\nicecal{F}_1$, $G\in\nicecal{G}_2$
where we only obtain $\approx_{\log\log}$ asymptotics.:

\begin{Th}\label{Th22.10a} \mbox{} \\
{\rm (1:1)} Assume $F\in\nicecal{F}_1,\,
G\in\nicecal{G}_1$. Then
$\displaystyle \overline H(x)\, \logsim\, \exp\process
{-c_{11}\log^{\theta_{11}}x}$ where
$\theta_{11}=\eta/\gamma$, $c_{11}=\alpha/\beta^{\theta_{11}}$;\\
{\rm (2:2)} Assume $F\in\nicecal{F}_2,\,
G\in\nicecal{G}_2$. Then
$\displaystyle \overline H(x)\, \logsim\, \frac{1}{x^{\theta_{22}}}
\,=\,\exp\process{-\theta_{22}\log x}$
where $\theta_{22}=\alpha/\beta$;\\
{\rm (2:1)} Assume $F\in\nicecal{F}_2,\,
G\in\nicecal{G}_1$. Then
$\displaystyle \overline H(x)\, \logsim\, 
\frac{1}{\log^{\theta_{21}}x}\,=\,\exp\process{-\theta_{21}\log\log x}$
where $\theta_{21}=\alpha/\gamma$;\\
{\rm (1:2)} Assume $F\in\nicecal{F}_1,\,
G\in\nicecal{G}_2$. Then
$\displaystyle \overline H(x)\, \approx_{\log \log}\, 
\exp\process{-x^{\theta_{12}}}$ where
$\theta_{12}=\eta/(\beta+\eta)\in(0,1)$.
\end{Th}
Note that the asymptotic expressions are in agreement with $\overline H(x)$
being necessarily heavy-tailed, cf.\ Proposition~\ref{Prop15.6a}.
E.g.\ the asymptotics in part (1:2) is as for the heavy-tailed
Weibull distribution, and the one in part (1:1) as for regular variation
if $\theta_{11}=1$ and as for the lognormal distribution if
$\theta_{11}=2$.

Generalizing \cite{Sheahan06}, we will also show:
\begin{Prop}\label{Prop28.5a}
Assume $g(t)\ge cf(t)^{1/\alpha-\epsilon} $ for all large $t$, where
$\alpha,\epsilon>0$, $c<\infty$. Then $\int_0^\infty x^\alpha H(\dd x)\ <\ 
\infty$. If $g(t)\le cf(t)^{1/\alpha} $ for all  large $t$, where
$\alpha>0$, $c>0$, then $\int_0^\infty x^\alpha H(\dd x)\ =\ 
\infty$. 
\end{Prop}
For example, the mean of $H$ is finite when the tail of $F$ is slightly
lighter than the tail of $G$ and infinite when it is equal or
or heavier. Similar, checking finite variance amounts to a comparison
of $\Fb$ and $\Gb^2$.

Our main results on sharp
asymptotics is as follows
(here and in the following, slowly varying functions are assumed
to have the additional property of being bounded on compact
subsets of $(0,\infty)$):
\begin{Th}\label{Prop12.3a}
Assume 
\begin{equation}\label{12.3a}
f(t)\ =\  g(t)\overline G(t)^{\beta-1}L_0\bigl(\overline G(t)\bigr)
\end{equation}
where $L_0(s)$ is slowly varying at $s=0$. Then
\begin{equation}\label{12.3b}
\overline H(x)\ \sim\ \frac{\Gamma(\beta)}{\mu^\beta}
\frac{L_0(1/x)}{ x^\beta},\ \
x\to\infty\,.
\end{equation}
\end{Th} 
Here $f(x)\sim g(x)$ means $f(x)/g(x)\to 1$. For example:
\begin{Cor}\label{Cor13.3a}
Assume $f,g$ belong to the class of regularly varying densities
of the form $L(t)/t^{1+\alpha}$ where $L$ is slowly varying,
with parameters $\alpha_F,L_F$ for $f$ and  $\alpha_G,L_G$ for $g$.
Then $\overline H(x)\ =\ L_H(x)/x^{\alpha_H}$, where $\alpha_H=
\alpha_F/\alpha_G$ and $L_H$ is slowly varying with
$$L_H(x)\ \sim\ \frac{\Gamma(\alpha_H)\alpha_G^{\alpha_H-1}}
{\mu^{\alpha_H}}
\frac{L_F\bigl(x^{1/\alpha_G}\bigr)}
{L_G^{\alpha_H}\bigl(x^{1/\alpha_G}\bigr)}\,.$$
\end{Cor}
\begin{Cor}\label{Cor13.3b}
Assume $f,g$ belong to the class of densities
of the form $\e^{-\lambda t^\eta}t^{\alpha}L(t)$ where $L$ is slowly varying
at $t=\infty$,
with parameters $\lambda_F,\alpha_F,L_F$ for $f$ and  
$\lambda_G,\alpha_G,L_G$ for $g$, and the same $\eta=\eta_F=\eta_G$.
Then $\overline H(x)\ =\ L_H(x)/x^{\alpha_H}$, where $\alpha_H=
\lambda_F/\lambda_G$ and $L_H$ is slowly varying with
$$L_H(x)\ \sim\ \frac{\Gamma(\alpha_H)\lambda_G^{\alpha_H-1-\omega}
\eta^{\alpha_H-1}}{\mu^{\alpha_H}
\lambda_G^{\alpha_F/\eta-
\alpha_G\alpha_H/\eta+\alpha_H-1}}
\log^\omega x\,\,
\frac{L_F\bigl(\log^{1/\eta}x\bigr)}
{L_G^{\alpha_H}\bigl(\log^{1/\eta}x\bigr)}\,,$$
where $\omega\,=\,\alpha_F/\eta+\alpha_H(\eta-\alpha_G-1)/\eta+1/\eta-1$.
\end{Cor}
Of course Corollary~\ref{Cor13.3a} is close in spirit to
Theorem~\ref{Th22.10a}(2:2); the conditions are slightly stronger,
but so are also the conclusions. The difference between
Corollary~\ref{Cor13.3b} and Theorem~\ref{Th22.10a}(1:1)
is somewhat more marked, since Corollary~\ref{Cor13.3b}
only applies when $\eta_F=\eta_G$ (i.e., $\eta=\gamma$ in the
notation of Theorem~\ref{Th22.10a}, where $\eta=\gamma$
is not required).

Finally consider ordering and comparison results.
One expects intuitively a heavier tail of  $F$ to lead to
a heavier $\overline H$. The precise
statement of this is in terms of stochastic order (s.o.):

\begin{Prop}\label{Prop17.6aw} 
Assume  given two task time distributions $F_1$, $F_2$ such
that $F_1$ is smaller than $F_2$ in s.o.,
that is, $\overline F_1(t)\le \overline F_2(t)$
for all $t$. Then also $H_{F_1,G}\le H_{F_2,G}$ in s.o.\ 
for any fixed $G$. 
\end{Prop}
This follows from \refs{16.5a} and the coupling characterization
of s.o.\ (\cite{Stoyan})
by noting that if $T_1\le T_2$, then (in obvious notation)
$N(T_1)\le N(T_2)$ and hence
$S(T_1)\le S(T_2)$, $X(T_1)\le X(T_2)$

Similarly, one expects a lighter tail of  $G$ to lead to
a larger $\tau$. However, stochastic ordering cannot be
inferred since if  $G_1$, $G_2$ are given ($F$ is fixed)
such that $G_1$ is smaller than $G_2$ in s.o,
then on one hand $N$ is smaller for $G_2$ than for $G_1$
for any $t$ but on the other the $U_i(t)$ are larger. 
However, we will establish an asymptotic order under 
a slightly stronger condition than $G_1$ being smaller than $G_2$ in s.o.:

\begin{Prop}\label{Prop17.6a} 
Assume  that $G_1$ is smaller than $G_2$ in s.o.\
and that in addition $\limsup_{t\to\infty}\overline G_1(t)/ \overline G_2(t)
<1$. 
Then for $F$ fixed,
$$
\limsup_{t\to\infty} \frac{\overline H_{F,G_2}(t)}{\overline H_{F,G_1}(t)}
\ \le\ 1\,.$$
\end{Prop}

\section{Geometric Sums. Bounded Job Time $T$} \label{S:Geom}

Given $T=t$, the number $N(t)$ of restarts is geometric with failure parameter
$G(t)=\Prob(U_i\le t)=1-\overline G(t)$ so that
$$\Prob\bigl(N(t)>n\bigr)\,=\,G(t)^n,\ \ 
\Exp N(t)\,=\,\frac{G(t)}{\overline G(t)}\,
\sim\, \frac{1}{\overline G(t)}\,.$$
It follows that given $T=t$, we can write
$$\tau\ \eqdistr\ t\,+\,S(t)\ \ \mbox{where}\ \ S(t)\,=\,
\sum_{i=1}^{N(t)}U_i(t)$$
(here $\eqdistr$ means equality in distribution)
where the $U_i(t)$ are independent of $N(t)$ and i.i.d.\ with the
distribution  $G_t$ being $G$ truncated to $[0,t)$, that is, with density
$G(t)^{-1}g(s)I(s\le t)$ at $s$. 
Then $\tau\eqdistr T+S(T)$ so that
\begin{equation}\label{12.6a}
\overline H(x)\ =\ \int_0^\infty \Prob\bigl(S(t)>x-t\bigr)\,f(t)\,
\dd t\,.
\end{equation}
This is the basic identity to be used in the following.

A first implication of \refs{12.6a} is that asymptotic
properties of geometric sums must play a role
for the asymtotics of $\overline H(x)$. We shall
use Cram\'er-Lundberg theory, cf.\ \cite{RP}, \cite{APQ}, \cite{Willmot}, 
more precisely the following result:

\begin{Prop}\label{Prop20.10a} 
Let $V_1,V_2,\ldots$ be i.i.d.\ with common density $k(v)$,
$N\in\nat$ an independent r.v.\ with $\Prob(N=n)=(1-\rho)\rho^n$,
and $S=V_1+\cdots+V_N$.
Then $\Prob(S>x)\sim C\e^{-\gamma x}$ where $\gamma$ is the solution
of $\rho\int_0^\infty \e^{\gamma y}k(y)\,\dd y\,=\,1$ and
$C=(1-\rho)/\gamma B$ where 
$B\,=\,\rho\int_0^\infty y\e^{\gamma y}k(y)\,\dd y$. Furthermore,
letting $$c_-(x)\,=\,\inf_{0\le z\le x, \overline K(z)>0}\frac{
e^{\gamma x}\overline K(x)}
{\int_x^\infty \e^{\gamma y}k(y)\,\dd y},\ \
c_+(x)\,=\,\sup_{0\le z\le x, \overline K(z)>0}\frac{
e^{\gamma x}\overline K(x)}
{\int_x^\infty \e^{\gamma y}k(y)\,\dd y,}$$
we have the Lundberg inequality
$$ c_-(x)\e^{-\gamma x}\ \le\ \Prob(S>x)\ \le\ 
c_+(x)\e^{-\gamma x}$$ for all $x$.
\end{Prop}
For a  proof, see Willmot \& Lin \cite{Willmot} pp.\ 108-109.
Alternatively, Proposition~\ref{Prop20.10a}
follows easily from 
\begin{eqnarray*}\Prob(S>x)&=&\Prob(N\ge 1,V_1>x)\,+\,
\int_0^x \Prob(S>x-y)\Prob(N\ge 1,V_1\in \dd y)\\&=&
\rho\overline K(x)\,+\,
\int_0^x \Prob(S>x-y)\rho k(y)\,\dd y\,,
\end{eqnarray*}
which is a defective renewal equation to which standard theory
applies (see \cite{APQ} V.7 and also \cite{RP} III.6c).
 \begin{Cor}\label{Cor20.10a} In the {\rm RESTART} setting,
$\Prob\bigl(S(t)>x\bigr)\,\sim\,C(t)\e^{-\gamma(t)x}$,
$x\to\infty$, where  $\gamma(t)>0$ is the solution of 
$\int_0^{t}\e^{\gamma(t) y} G(\dd y)\,=\,1$ 
and $C(t)\,=\,\overline G(t)/\gamma(t)B(t)$
where  $B(t)\,=\,\int_0^{t}y\e^{\gamma(t) y}g(y)\,\dd y$.
This estimate is uniform in $t_1\le t\le t_2$ for given
$0<t_1<t_2$.
Furthermore,
$$\e^{-\gamma(t) t}\e^{-\gamma(t)x}\ \le\ \Prob\bigl(S(t)>x\bigr)\ \le\ 
\e^{-\gamma(t) t}\,.
$$
\end{Cor} 
\proof The first statement is a trivial translation of
the first statement of Proposition~\ref{Prop20.10a}. For the two-sided
Lundberg inequality, note that 
in the {\rm RESTART} setting with $K(y)=K_t(y)=\Prob\bigl(U(t)\le y\bigr)$, 
the integral in the definition
of $c_-(x)$ extends only up to $t$ which gives $c_-(x)\ge \e^{-\gamma t}$,
and that $c_+(x)\le 1$.  For the uniformity of the Cram\'er-Lundberg
approximation, appeal to uniform estimates of the
renewal functions corresponding to the $\e^{\gamma(t) y}K_t(\dd y)$
as given, e.g., Kartashov~\cite{Karta0}, \cite{Karta}
(see also Wang \& Woodroofe~\cite{WW}).\halmos

\smallskip

In particular, Corollary~\ref{Cor20.10a} settles
the case of a fixed job size:
\begin{Cor}\label{Cor16.5a} Assume $T\equiv t_0$ and $\Gb(t_0)>0$. Then
$$\overline H(x)\ \sim\ C(t_0)\e^{\gamma(t_0)t_0}\e^{-\gamma(t_0)x}
\,.$$
\end{Cor}

In the case of an infinite support of $f$, Corollary~\ref{Cor16.5a}
shows that the tail of $H$ is  heavier than $\e^{-\gamma(t)x}$ for all
$t$ (note that $\gamma(t)\downarrow 0$ as $t\to\infty$; more precise
estimates are given later). This observation proves
Proposition \ref{Prop15.6a}.

If $T$ is random, we need to mix over $t$
with weights $f(t)$. If the support of $f$ has a finite upper
endpoint $t_0$, Corollary~\ref{Cor16.5a} suggests that
the asymptotics of $\overline H(x)$ is not too far from
$\e^{-\gamma(t_0)x}$, and in fact, we shall show:

\begin{Prop}\label{Prop12.10b}
Assume that the support of $F$ has upper endpoint $0<t_0<\infty$,
that $\overline G(t_0)>0$
and that
\begin{equation}\label{12.10c} f(t)\ \sim\ A(t_0-t)^\alpha,\ \ t\uparrow t_0,
\end{equation}
for some $0<A<\infty$ and some $\alpha\ge 0$. Then $$\overline H(x)\ \sim\
\frac{AB(t_0)^\alpha\overline G(t_0)\Gamma(\alpha+1)}
{\gamma(t_0)\e^{\alpha\gamma(t_0)}g(t_0)^{\alpha+1}}
\frac{\e^{-\gamma(t_0) x}}{x^{\alpha+1}}\,.$$ 
\end{Prop}
\proof
For simplicity of notation, write $B=B(t_0)$, $\gamma=\gamma(t_0)$ etc.

It is easy to see that $\gamma(t)$ is continuous and differentiable in
$t$. To obtain the asymptotics as $t\uparrow t_0$ we write
\begin{eqnarray*} 1&=&\int_0^{t}\e^{\gamma(t) y}G(\dd y)\ =\
\int_0^{t_0}\e^{\gamma(t) y}G(\dd y)\,-\,
\int_t^{t_0} \e^{\gamma(t) y}G(\dd y)\\
&&\int_0^{t_0}\e^{\gamma y}\bigl[1+(\gamma(t)-\gamma)y\bigr]\,G(\dd y)\,-\,
(t_0-t)\e^{\gamma t_0}g(t_0)\,+\,\oh\bigl(\gamma(t)-\gamma\bigr)\\
&=& 1\,+\, (\gamma(t)-\gamma)B\,-\,
(t_0-t)\e^{\gamma t_0}g(t_0)\,+\,\oh\bigl(\gamma(t)-\gamma\bigr)\\
\end{eqnarray*}so that 
\begin{equation}\label{18.10a}
\gamma(t)-\gamma\,\sim\,(t_0-t)D
\end{equation}
where $D=\e^{\gamma t_0}g(t_0)/B$.
Appealing to the uniformity in Corollary~\ref{Cor20.10a}, we therefore get
\begin{eqnarray*}\overline H(x)&=& \int_0^{t_0}Z_t(x-t)f(t)\,\dd t\ =\
\int_{t_0-\epsilon}^{t_0}Z_t(x-t)f(t)\,\dd t \ +\ \oh(\e^{-\gamma x})\\
&=&r_1(\epsilon)\int_{t_0-\epsilon}^{t_0}C
\e^{-\gamma(t)(x-t)}A(t_0-t)^\alpha\,\dd t
 \ +\ \oh(\e^{-\gamma x})
\\ &=&r_2(\epsilon)AC\e^{-\gamma (x-t_0)}
\int_{t_0-\epsilon}^{t_0}\e^{-(\gamma(t)-\gamma)x}(t_0-t)^\alpha\,\dd t
\ +\ \oh(\e^{-\gamma x})
\end{eqnarray*}
where $Z_t(x)=P(S(t)>x)$ and $r_1(\epsilon),r_2(\epsilon),\ldots 
\to 1$ as $\epsilon\downarrow 0$. Thus
substituting $y=(\gamma(t)-\gamma)x$ and noting that
$\dd y\sim -D\dd t$ by \refs{18.10a}, we get up 
to the $\oh(\e^{-\gamma x})$ term that
\begin{eqnarray*}\overline H(x)&=&
r_2(\epsilon)ACD^{-\alpha-1}\frac{\e^{-\gamma (x-t_0)}}{x^{\alpha+1}}
\int_0^{(\gamma(t_0-\epsilon)-\gamma)x}y^\alpha \e^{-y}
\,\dd y
\end{eqnarray*}
Letting first $x\to\infty$, next $\epsilon\downarrow 0$,
and rewriting the constants
completes the proof.\halmos

\section{Proof of Lemma~\ref{L20.10c}}\label{S:First}

We will need the asymptotics of the Cram\'er root $\gamma(t)$:
\begin{Lemma}\label{L20.10b} As $t\to\infty$,
$\gamma(t)\,\sim\,\mu\overline G(t)$.
\end{Lemma}
\proof Consider
\begin{equation}\label{bed31.5a}
\int_0^t\bigl(\e^{\gamma(t) y}-1-\gamma(t)y\bigr)\,G(\dd y)\ =\ 
\Gb(t)-\gamma(t)\bigl(1/\mu-\oh(1)\bigr)\,.
\end{equation}
The non-negativity of the l.h.s.\ yields $\gamma(t)=\Oh\bigl(\Gb(t)\bigr)$.
Since $t\Gb(t)$ because of $\mu>0$, the integrand in \refs{bed31.5a}
can therefore be writtes as $\gamma(t)y\epsilon(y,t)$ where $\epsilon(y,t)
\to 0$ uniformly in $y\le t$ as $t\to\infty$. Therefore \refs{bed31.5a}
equals $\gamma(t)\oh(1)$ which shows the assertion.\halmos

\smallskip

\noindent {\em Proof of Proposition }\ref{Prop15.6a}.
Given $\epsilon>0$, choose $t_0$ such that $\gamma(t_0)<\epsilon$,
cf.\ Lemma~\ref{L20.10b}, and $a$ so large that
$\gamma(t_0+a)<\gamma(t_0)$. We then get
\begin{eqnarray*}\lefteqn{\liminf_{x\to\infty}\e^{\epsilon x}\overline H(x)
\ \ge \ \liminf_{x\to\infty}\e^{\gamma(t_0) x}\overline H(x)\ \ge\
\liminf_{x\to\infty}\int_{t_0+a}^\infty 
\frac{\Prob\bigl(S(t)>x\bigr)}{\e^{-\gamma(t_0) x}}}\\
&\ge&\int_{t_0+a}^\infty \liminf_{x\to\infty}
\frac{\Prob\bigl(S(t)>x\bigr)}{\e^{-\gamma(t_0) x}}f(t)\,\dd t\ =\ 
\int_{t_0+a}^\infty\infty\,\cdot\,f(t)\,\dd t\,=\,\infty
\end{eqnarray*}
where we used Fatou's lemma in the third step and
Corollary~\ref{Cor20.10a} in the next.~\halmos

\begin{Lemma}\label{L7.6a} For any $t_0<\infty$,
$\Prob\bigl(\tau>x,\,T\le t_0\bigr)$ goes to zero
at least exponentially fast.
\end{Lemma}
\proof By Lundberg's inequality,
$$\Prob\bigl(\tau>x,\,T\le t_0\bigr)\ \le\ 
F(t_0)\Prob\bigl(S(t_0)>x-t_0\bigr)\ \le\ 
F(t_0)\e^{-\gamma(t_0)(x-t_0)}\,.\eqno{\halmos}$$

\begin{Lemma}\label{L28.5a}
Define $S_n(t)=U_1(t)+\cdots+U_n(t)$, 
$m(t)=\Exp U_i(t)=\Exp[U_i\,|\,U_i\le t]$. Then
$$\Prob\Bigl(\bigl|S_n(t)/n-m(t)\bigr|>\epsilon
\Bigr)\ =\ \oh(1),\ \ n\to\infty,$$
where the $\oh(1)$ is uniform in $t>\delta$ for any $\delta>0$. 
\end{Lemma}
\proof Define $U_i(t,n)=U_i(t)I\bigl(U_i(t)<n\bigr)$. Then, in obvious notation
$$\Prob\bigl(S_n(t,n)\ne S_n(t)\bigr)\ \le\ n
\Prob\bigl(U_n(t,n)\ne U_n(t)\bigr)\ \le\
\frac{n}{G(t)}\Prob(U>n)
$$  
goes to zero uniformly in  $t>\delta$ because of $\Exp U<\infty$.
Further,
$$\frac{1}{n}\Exp U_i(t,n)^2 \ =\ \int_0^n \frac{2x}{n}
\Prob\bigl(U_i(t,n)>x\bigr)
\, \dd x \le\ \frac{2}{G(t)}\int_0^n \frac{x}{n}
\Prob\bigl(U>x\bigr)\,\dd x \ =\ \oh(1)
$$
uniformly in  $t>\delta$, as follows by dominated convergence
with  $\Prob\bigl(U>x\bigr)$ as majorant. Hence by Chebycheff's
inequality,
$$\Prob\Bigl(\bigl|S_n(t,n)/n-m(t,n)\bigr|>\epsilon
\Bigr)\ \le\ \frac{n\Exp U_i(t,n)^2}{n^2\epsilon^2}\ =\ \oh(1)
$$
uniformly in  $t>\delta$. Also
$$m(t)-m(t,n)\ = \ \Exp U_i(t)I\bigl(U_i(t)\ge n\bigr)\ \le\ 
\frac{1}{G(t)}\Exp U I\bigl(U\ge n\bigr)\ =\ \oh(1)
$$
uniformly in  $t>\delta$. Putting these estimates together
completes the proof.
\halmos

\smallskip

\noindent{\em Proof of Lemma }\ref{L20.10c}. Given $\epsilon>0$,
it follows by Lemma~\ref{L20.10b} that we can
choose $t_0$ such that $\gamma(t)\ge\mu\overline G(t)(1-\epsilon)$
and (since $G$ has finite mean)
$\gamma(t)t<\log(1+\epsilon)$ for $t\ge t_0$ . 
Thus 
by the upper Lundberg bound and Lemma~\ref{L7.6a},
\begin{eqnarray*}
\overline H(x)&=& \int_{t_0}^\infty \Prob\bigl(S(t)>x-t\bigr)
f(t)\,\dd t\ +\ \oh(\e^{-rx})\\&\le& \int_{t_0}^\infty \e^{-\gamma(t)(x-t)}
f(t)\,\dd t\ +\ \oh(\e^{-rx})\\ &\le& (1+\epsilon)\int_{t_0}^\infty 
\e^{-\mu\overline G(t)(1-\epsilon)x}f(t)\,\dd t\ +\ \oh(\e^{-rx})\\& \le&
(1+\epsilon)I_-(x,\epsilon)\ +\ \oh(\e^{-rx})
\end{eqnarray*}
for some $r>0$.
Now note that $\overline H(x)$ decays slower than $\e^{-rx}$
by  Proposition \ref{Prop15.6a}.

For the lower bound, let $\epsilon>0$ be given and let 
$\epsilon_1,\epsilon_2>0$ satisfy $(1+\epsilon)(1-\epsilon_1)$
$>1+\epsilon_2>1$.
Lemma~\ref{L28.5a} implies that  there is an $n_0$ such that
$$\Prob\bigl(S_n(t)>nm(t)(1-2\epsilon_1)\bigr)\ \ge\ 1-\epsilon
$$
for all $n\ge n_0$ and all $t\ge t_0$.
Since $m(t)\to 1/\mu$, we have then also
$$\Prob\bigl(S_n(t)>n(1-\epsilon_1)/\mu\bigr)\ \ge\ 1-\epsilon $$
for all $n\ge n_0$ and all $t\ge t_0$. Choose next $\overline g_0$
such that $\e^{-(1+\epsilon_2)\overline g}$ $\le 1-\overline g$
for $0<\overline g<\overline g_0$. Replacing $t_0$ by a larger $t_0$
if necessary, we may assume $\overline G(t)<\overline g_0$
for $t\ge t_0$ and get
\begin{eqnarray*}
\overline H(x)&\ge& \int_{t_0}^\infty \Prob\bigl(S(t)>x\bigr)
f(t)\,\dd t\\ &\ge& (1-\epsilon)\int_{t_0}^\infty \Prob\bigl(N(t)>x\mu/
(1-\epsilon_1)\bigr) f(t)\,\dd t\\
&=&(1-\epsilon)\int_{t_0}^\infty G(t)^{x\mu/(1-\epsilon_1)} f(t)\,\dd t
\\ & \ge& (1-\epsilon)\int_{t_0}^\infty \exp\process{-\overline G(t)x\mu
(1+\epsilon_2)/(1-\epsilon_1)} f(t)\,\dd t\\
&\ge&(1-\epsilon)\int_{t_0}^\infty \exp\proces{-\Gb(t)x\mu
(1+\epsilon)} f(t)\,\dd t\,.
\end{eqnarray*}
Since the last integral differs from $I_+(x,\epsilon)$ by a term
which goes to zero exponentially fast and hence is $\oh(\overline H(x))$,
the proof is complete.\halmos

\smallskip

\noindent {\em Proof of Proposition} \ref{Prop02.01a}.
When $F=G$, we have $\dd\overline G(t)=-f(t)$. Hence
\begin{eqnarray*}I_\pm(x,\epsilon) &=& \Bigl[\frac{1}{\mu x(1\pm\epsilon)}
\e^{-\mu\overline G(t)x(1\pm\epsilon)}\Bigr]_{t_0}^\infty\ =\ 
\frac{1}{\mu x(1\pm\epsilon)}
\bigl(1-\e^{-\mu\overline G(t_0)x(1\pm\epsilon)}\bigr)\\ & \sim& 
\frac{1}{\mu x(1\pm\epsilon)}\,.
\end{eqnarray*}
 The assertion now follows easily from Lemma~\ref{L20.10c} by letting
first $x\to\infty$ and next $\epsilon\to 0$.\halmos

\smallskip

\noindent {\em Proof of Proposition} \ref{Prop28.5a}.
Under the assumptions of the last part of the Proposition,
$\Gb(t)\le c_1\Fb(t)^{1/\alpha}$ for $t\ge t_0$ and hence
\begin{eqnarray*}
\int_0^\infty x^\alpha H(\dd x)&=&\alpha\int_0^\infty x^{\alpha-1}
\Hb(x)\,\dd x\\
&\ge&c_2\int_0^\infty x^{\alpha-1}\dd x\int_0^\infty \e^{-\mu\Gb(t)x}
f(t)\,\dd t\\
&\ge&c_3\int_{t_0}^\infty\frac{1}{\Gb(t)^\alpha}f(t)\,\dd t\ \le\ 
c_4\int_{t_0}^\infty\frac{1}{\Fb(t)}f(t)\,\dd t\\
&=&c_4\int_0^1 \frac{1}{y}\dd y\ =\ \infty\,,
\end{eqnarray*}
proving the last part of the Proposition. For the first part, we get
similarly
$$\int_0^\infty x^\alpha H(\dd x)\ \le\ 
c_5\,+\,c_6\int_0^1 \frac{1}{y^{1-\alpha\epsilon}}\dd y\ <\ \infty\,.
\eqno{\halmos}$$

\smallskip

\noindent {\em Proof of Proposition} \ref{Prop17.6a}. 
Write $I^1_\pm(x,\epsilon)$, $\mu_1$ when $G=G_1$ and similarly
for $G_2$. We may assume $G_1\ne G_2$. Then the s.o.\ assumption
implies $\mu_1<\mu_2$. Hence if
$\epsilon>0$ is so small that $\mu_1(1+\epsilon)<\mu_2(1-\epsilon)$, 
we have 
$\mu_1\overline G_1(t)(1+\epsilon)$ $\le$ $\mu_1\overline G_2(t)(1-\epsilon)$
for all $t$. We then obtain
$$\overline H_{G_1}(x)\ \ge\ (1-\epsilon)I^1_+(x,\epsilon)
\ \ge\ (1-\epsilon)I^2_-(x,\epsilon)\ \ge\ 
\frac{1-\epsilon}{1+\epsilon}\overline H_{G_2}(x)\,,$$
where the outer inequalities are asymptotic and the inner one exact.
Let first $x\to\infty$ and next $\epsilon\to 0$.\halmos

\section{Proofs: Integral Asymptotics}\label{S:Further}

\begin{Lemma}\label{L21.10a}
For given constants $a,b,\gamma,\eta>0$,
$$I\ =\ \int_{t_0}^\infty \exp\process{-\e^{-bt^\gamma}z-at^\eta}\,\dd t
\ \logsim\ \e^{-ab^{-\eta/\gamma}\log^{\eta/\gamma}z}$$
as $z\to\infty$.
\end{Lemma}
\proof Let $c=ab^{-\eta/\gamma}$, $t_1=(\log z/b)^{1/\gamma}$ and let $I_1,I_2$
be the contributions to $I$ from the intervals $(t_1,\infty)$,
resp.\ $(t_0,t_1)$. In $I_1$, we bound the first term 
in the exponent below by $0$ so 
$$I_1\ \le\ \int_{t_1}^\infty 
\e^{-at^\eta}\,\dd t\ \logsim\ \e^{-at_1^\eta}\ =\ 
\e^{-c\log^{\eta/\gamma}z}\, .$$
In $I_2$, we substitute $y=\e^{-bt^\gamma}z$. Then
$$t\,=\,(\log z-\log y)^{1/\gamma}b^{-1/\gamma},\ \
\dd t\,=\,-\frac{1}{\gamma y}(\log z-\log y)^{1/\gamma-1}b^{-1/\gamma}\dd y$$
so that $I_2$ becomes
\begin{eqnarray*}&&\frac{1}{\gamma b^{1/\gamma}}\int_1^{\e^{-t_0^\gamma}z}\frac{1}{y}
(\log z-\log y)^{1/\gamma-1}
\exp\process{-y-c
(\log z-\log y)^{\eta/\gamma}}\,\dd y
\end{eqnarray*}
We split this integral into the contributions $I_3,I_4$
from the intervals $[1,2)$, $[2,\e^{-t_0^\gamma}z)$.
Here 
$$I_3 \sim\ \frac{1}{\gamma b^{1/\gamma}}
\log^{1/\gamma-1}z\,\,\e^{-c\log^{\eta/\gamma}z} \int_1^2
\frac{1}{y}\e^{-y}\,\dd y\ \logsim \ \e^{-c\log^{\eta/\gamma}z}\,.
$$
For $I_4$, we write $q=\eta/\gamma$, $h(y)=-y/2+c\log^qz-c(\log z-
\log y)^q$.
Then 
$$h(y)\ =\ -y/2+cq\int_1^y \frac{(\log z-
\log y)^{q-1}}{y}\,\dd y\ \le\ -y/2+cq\log^{q-1} z\log y\,.
$$
The r.h.s.\ is maximized for $y_z\,=\,2cq\log^{q-1} z$, where
$$h(y_z)\ =\ \log^{q-1} z\Oh(\log\log z)\ =\ \oh\bigl(\log^q z\bigr)\,.$$
Hence
\begin{eqnarray*}I_4&\le&b^{-1/\gamma}\log^{1/\gamma-1}z\e^{-c\log^q z}
\int_1^z\frac{1}{y}\e^{-y/2+h(y)}\,\dd y\\
&\le&\e^{-\bigl(c+\oh(1)\bigr)\log^q z} \int_1^z\frac{1}{y}\e^{-y/2}\,\dd y
\ \logsim\ \e^{-c\log^q z}\,.
\end{eqnarray*}

Adding these estimates shows that $ \e^{-ab^{-\eta/\gamma}\log^{\eta/\gamma}z}$
is an asymptotic upper bound in the logarithmic sense, and that it is also
a lower one follows from the estimate for $I_3$.
\halmos 

\begin{Lemma}\label{L26.10a}
For given constants $a,b$,
$$\int_{t_0}^\infty \e^{-t^{-\beta} z}\frac{1}{t^{\alpha+1}}\,\dd t
\ \logsim\ \frac{1}{z^{\alpha/\beta}}\,.
$$
\end{Lemma}
\proof Substitute $y=t^{-\beta} z$ to get
$$\int_{t_0}^\infty 
 \e^{-t^{-\beta} z}\frac{1}{t^{\alpha+1}}\,\dd t\ =\ 
\int_0^{t_0^{-\beta} z}
\frac{y^{\alpha/\beta-1}}{\beta z^{\alpha/\beta}}
\e^{-y}\,\dd y\ \sim \
\frac{\Gamma\bigl(\alpha/\beta\bigr)}{\beta z^{\alpha/\beta}}
)\ \logsim\ \frac{1}{z^{\alpha/\beta}}\,.
\eqno{\halmos}$$

\begin{Lemma}\label{L21.10aa}
For given constants $a,b,\gamma>0$,
$$I\ =\ \int_{t_0}^\infty \exp\process{-\e^{-bt^\gamma}z-(a+1)\log t}\,\dd t
\ \logsim\ \frac{1}{\log^{a/\gamma} z}$$
as $z\to\infty$.
\end{Lemma}
\proof Let again $t_1=t_1(z)=\bigl(\log z/b\bigr)^{1/\gamma}$ 
(then $\e^{-bt_1^\gamma}z=1$)
and let $I_1,I_2$
be the contributions to $I$ from the intervals $(t_1,\infty)$,
resp.\ $(t_0,t_1)$. In $I_1$, $0\,\le\,\e^{-bt^\gamma}z\,\le 1$
and so
$$\frac{\e^{-1}}{a\bigl(\log(z/b)\bigr)^{a/\gamma} } \ \le\ 
I_1\ \le\ \frac{1}{a\bigl(\log(z/b)\bigr)^{a/\gamma} }\,.$$

For $I_2$, let $r_z(t)=\e^{-bt^\gamma}z+(a+1)\log t$, 
$s(t)=t^\gamma \e^{-bt^\gamma}$. Then
$$r'_z(t)\ =\ -b\gamma t^{\gamma -1}\e^{-bt^\gamma}z\,+\,
\frac{a+1}{t}\ =\ \frac{1}{t}\bigl[-b\gamma s(t)z+(a+1)\bigr]\,.$$
Since $s$ is continuous with $s(t_0)>0$ and $s(t)$ is monotonically decreasing 
for large $t$ with limit 0, 
we have $s(t)\ge s(t_1)=\log (z/b)/z$ 
for all $t_0\le t\le t_1$ and all large $z$
because of $t_1(z)\to\infty$. Hence $r'_z(t)<0$ for 
$t_0\le t\le t_1$ and all large $z$ so that
\begin{eqnarray*}I_2& =& \int_{t_0}^{t_1} \e^{-r_z(t)}\,\dd t\ \le\ 
(t_1-t_0)\e^{-r_z(t_1)}\\
&\le& t_1\e^{-(a+1)\log t_1}\ =\ 
\frac{1}{t_1^a}\ =\ \frac{1}{\bigl(\log(z/b)\bigr)^{a/\gamma}}\,.
\end{eqnarray*}
Putting the upper bounds for $I_1,I_2$ together and noting that
$\bigl(\log(z/b)\bigr)^{a/\gamma}\,\logsim\, \log^{a/\gamma}x$
shows that $\log^{-a/\gamma}z$ is an upper bound in the
logarithmic sense, and that it is also a lower bound follows from
the lower bound for $I_1$.
\halmos

\begin{Lemma}\label{L21.10aaa}
Let $\eta>0$ be fixed. Then for any  $a,b>0$,
$$I\ =\ \int_{t_0}^\infty \exp\process{-t^{-b}z-at^\eta}\,\dd t
\ \logsim\ \exp\process{-c_{12}(a,b)z^{\theta_{12}(b)}}$$
as $z\to\infty$ where $$\theta_{12}(b)=\eta/(b+\eta),\ \  
c_{12}(a,b)\ =\ a^{1-\theta_{12}}\bigl[
(\eta/b)^{1-\theta_{12}(b)}+(b/\eta)^{\theta_{12}}\bigl]\,.$$
\end{Lemma}
\proof We choose $t_1=t_1(z)$ to minimize $f(t)\,=\,t^{-b}z+at^\eta$
which gives
$$t_1\ =\ \Bigl(\frac{b z}{a\eta}\Bigr)^{1/(b+\eta)}
\,, \ \ f(t_1)\ =\ c_{12}(a,b)z^{\theta_{12}(b)}\,.$$
Thus the claim of the lemma can be written as $I\,\logsim\,\e^{-f(t_1)}$.
As lower bound, we use
$$\int_{t_1}^{t_1+1} \exp\process{-t^{-b}z-at^\eta}\,\dd t\ \ge\
\exp\process{-t_1^{-b}z-a(t_1+1)^\eta}\ \logsim\ 
\e^{-f(t_1)}$$
where in the last step we used $(t+1)^\eta=t^\eta\bigl(1+\oh(1)\bigr)$.
For the upper bound, we write $I=I_1+I_2+I_3$ where $I_1,I_2,I_3$
are the contributions from the intervals $t_0<t<t_1$, $t_1<t<Kt_1$, resp.\
$Kt_1<t<\infty$ where $K$ satisfies $aK^\eta>c_{12}(a,b)$. 
Since $f$ is decreasing in the interval $t_0<t<t_1$
and increasing in $t_1<t<\infty$, 
we have $I_1\le t_1\e^{-f(t_1)}\logsim f(t_1)$ and 
$I_2\le (K-1)t_1\e^{-f(t_1)}\logsim \e^{-f(t_1)}$. 
Finally,
$$I_3\ \le\ \int_{Kt_1}^\infty \e^{-at^\eta}\,\dd t\ \logsim\
\e^{-aK^\eta t_1^\eta}$$
can be neglected because of the choice of $K$.\halmos

\smallskip

\noindent {\em Proof of Theorem} \ref{Th22.10a}.
In (1:1), we can choose $t_0$ such that $\overline G(t)\le\e^{-bt^\gamma}$
and $f(t)\ge\e^{-at^\eta}$, $t\ge t_0$, 
for any given $b<\beta$ and $a>\alpha$. With $I_\pm(x,\epsilon)$
as in Lemma~\ref{L20.10c}, we then get
\begin{eqnarray*}\lefteqn{\liminf_{x\to\infty}
\frac{\log I_+(x,\epsilon)}{-\log^{\eta/\gamma}x}}\\
&=&\liminf_{x\to\infty}
\frac{1}{-\log^{\eta/\gamma}x}\log\int_{t_0}^\infty \exp\process{
-\mu\overline G(t)x(1+\epsilon)}f(t)\, \dd t\\ &\ge&
\liminf_{x\to\infty}
\frac{1}{-\log^{\eta/\gamma}x}\log\int_{t_0}^\infty \exp\process{
-\mu \e^{-bt^{\gamma}}x(1+\epsilon)-at^\eta}\, \dd t\\ &=&
\liminf_{x\to\infty}
\frac{ab^{-\eta/\gamma}\log^{\eta/\gamma}\bigl(\mu x(1+\epsilon)\bigr)}
{-\log^{\eta/\gamma}x}\ =\ ab^{-\eta/\gamma}
\end{eqnarray*}
where we used Lemma~\ref{L21.10a} with $z=\mu x(1+\epsilon) $
in the third step. Letting $a\downarrow \alpha$, $b\uparrow \alpha$
shows that $\e^{-c_{11}\log^{\theta_{11}}x}$ is an asymptotic lower bound
in the logarithmic sense. That it is also an asymptotic upper bound
follows in the same way by noting that the contribution to $\overline H(x)$
from $(0,t_0)$ goes to zero exponentially fast by
Proposition~\ref{Prop15.6a}  for any $t_0$ and hence is negligible
compared to  $\e^{-c_{11}\log^{\theta_{11}}x}$.

Parts (2:2) and (2:1) follow in a similar way from 
Lemmas~\ref{L26.10a} and \ref{L21.10aa}. For (1:2), we choose
$b>\beta$, $a<\alpha$ and get
\begin{eqnarray*}\lefteqn{\liminf_{x\to\infty}
\frac{\log\bigl(-\log \overline H(x)\bigr)}
{\log\bigl(-\log\e^{-x^{\theta_{12}}}\bigr)}\ \ge \
\liminf_{x\to\infty}
\frac{\log\bigl(-\log I_+(x,\epsilon)\bigr)}{\theta_{12}\log x}}\\
&=&\liminf_{x\to\infty}
\frac{1}{\theta_{12}\log x}\Bigl(-\log\int_{t_0}^\infty \exp\process{
-\mu\overline G(t)x(1+\epsilon)}f(t)\, \dd t\Bigr)\\ &\ge&
\liminf_{x\to\infty}
\frac{1}{\theta_{12}\log x}\Bigl(-\log\int_{t_0}^\infty \exp\process{
-\mu t^{-b}x(1+\epsilon)-at^\eta}\, \dd t\Bigr)\\ &=&
\liminf_{x\to\infty}
\frac{\log\bigl(c_{12}(a,b)
\bigl(\mu x(1+\epsilon)\bigr)^{\theta_{12}(b)}\bigr)}
{\theta_{12}\log x}
\ =\ \frac{\theta_{12}(b)}{\theta_{12}}\,.
\end{eqnarray*}
Letting $a\uparrow \alpha$, $b\downarrow \alpha$
shows that $\e^{-x^{\theta_{12}}}$ is an asymptotic lower bound
in the $\approx_{\log\log}$ sense. That it is also an asymptotic upper bound
follows similarly.
\halmos

\noindent
{\em Proof of Theorem} \ref{Prop12.3a}.
In Lemma~\ref{L20.10c}, we insert \refs{12.3a}
and substitute $s=\Gb(t)$ to get
\begin{eqnarray*}I_\pm&=&\int_0^{s_0}\exp\proces{-sx\mu(1\pm\epsilon)}
s^{\beta-1}L_0(s)\,\dd s\,.
\end{eqnarray*}
where $s_0=\Gb^{-1}(t_0)$. Then by Karamata's Tauberian theorem
(\cite[Theorems 1.5.11 and 1.7.1]{BGT}),
$$I_\pm\ \sim\ \Gamma(\beta)
\frac{L\bigl(1/(x\mu(1\pm\epsilon)\bigr)}{x^\beta\mu^\beta
(1\pm\epsilon)^\beta}\ \sim\ 
\Gamma(\beta)\frac{L(1/x)}{x^\beta\mu^\beta
(1\pm\epsilon)^\beta}\,.
$$
Let $\epsilon\downarrow 0$.
\halmos

\bigskip \noindent
{\em Proof of Corollary} \ref{Cor13.3a}.
We have $\Gb(t)=L'_G(t)/t^{\alpha_G}$,
where $L'_G(t)\sim L_G(t)/(\alpha_G+1)$ as $t\to\infty$. Then
\refs{12.3a} holds with $\beta=\alpha_F/\alpha_G$ and
$L_0\bigl(\Gb(t)\bigr)=L_F(t)/L_G(t){L'_G{}}^{\beta-1}(t)$.
Note that $L_0$ is s.v.\ because the inverse of a
s.v.\ function is again s.v.\ (\cite[p.\ 28]{BGT}) and because the
composition of two s.v.\ functions is again s.v. Further,
(\cite[p.\ 29]{BGT})
$L_F\bigl(\Gb^{-1}(s)\bigr)\sim L_F\bigl(s^{-1/\alpha_G}\bigr)$ as
$s\downarrow 0$
and similarly for $L_G$. Thus, $L_0(s)\sim L_F\bigl(s^{-1/\alpha_G}\bigr)
 (\alpha_G+1)^{\beta-1}/L_G^\beta\bigl(s^{-1/\alpha_G}\bigr)$.
\halmos

\bigskip \noindent
{\em Proof of Corollary} \ref{Cor13.3b}.
We have
$$\Gb(t)\ \ =\ \int_t^\infty \e^{-\lambda_G y^\eta}y^{\alpha_G}L_G(y)\,\dd y
\sim\ \frac{1}{\eta\lambda_G}
\e^{-\lambda_G t^\eta}t^{\alpha_G+1-\eta_G}L_G(t)$$
(e.g., substitute $z=\e^{-\lambda_G y^\eta}$ and apply
Karamata's theorem).
>From this it is easy to see that
$$\Gb^{-1}(s)\ \sim\ (-\log s)^{1/\eta}/\lambda_G^{1/\eta}\,,\ \
s\downarrow 0.$$
In particular, $L_F\bigl(\Gb^{-1}(s)\bigr)\,\sim\,L_F\bigl((-\log s)^{1/\eta}
\bigr)$
and similarly for $L_G$.
Thus if $L_0$ is defined by \refs{12.3a} with 
$\beta=\lambda_F/\lambda_G$, we have
\begin{eqnarray*}L_0(s)&\sim&
\lambda_G^{\beta-1-\omega}\eta^{\beta-1}(-\log s)^\omega
\frac{L_F\bigl((-\log s)^{1/\eta}\bigr)}{L_G^\beta
\bigl((-\log s)^{1/\eta}\bigr)}
\end{eqnarray*}
which in particular shows that $L_0$ is s.v.\ at $s=0$.
Now just replace $s$ by $1/x$ and $\beta$ by
$\alpha_H$ to obtain the Corollary.
\halmos

\smallskip

\noindent
Corollary~\ref{Cor23.5a} is a special case of Corollary \ref{Cor13.3b}.

\section{Concluding Remarks} \label{S:Concl}

\noindent {\bf A)} The representation \refs{16.5a}
easily gives a proof of the asymptotics $\Hb(x)\sim 1/\mu x$
for the diagonal case $F=G$ (Proposition~\ref{Prop02.01a}). Indeed,
the event $N=n$ corresponds to the ordering $U_1<T,\ldots,U_n<T$,
$U_{n+1}>T$. Since the $n+2$ random variables $T,U_1,\ldots,U_n,
U_{n+1}$ are i.i.d.\ when $F=G$, we therefore have $\Prob(N=n)=
1/(n+2)(n+1)$, $\Prob(N>n)=1/(n+2)$
(we are grateful to Clive Anderson for a remark triggering this observation). 
One can now argue that
in order for $X$ to be large, $N$ has to be large which in turn
is only possible if $T$ is large. Then the distribution of the
$U_i$ is close to $G$, so that the geometric sum is approximately
$N/\mu$. Since $\Exp T=1/\mu<\infty$ implies that the tail
of $T$ is lighter than $1/x$, we therefore get
$$\Prob(X>x)\ \approx\ \Prob\Bigl(\sum_{i=1}^{N}U_i>xBbigr)\
\approx\ \Prob(N>\mu x)\ \sim\ \frac{1}{\mu x}\,.$$

Th argument is not hard to make rigorous, but we omit the details
since the further results of the paper are much more general than
Proposition~\ref{Prop02.01a} and require different proofs.

\medskip

\noindent {\bf B)}
An application of the above results occurs in
parallel computing. Assume that a job of length $NT$ is split into $N$ subjobs
of length $T$
which are placed on $N$ parallel processors ($N$ may run in the order of
hundreds or thousands). If one processor
fails, the corresponding subjob is restarted on a new processor.
With $X_1,\ldots,X_N$ the total times of the subjobs, the total
job time is then $M_N\,=\,\max(X_1,\ldots,X_N)$. The asymptotic
behaviour of $M_N$ as $N\to\infty$
is available from extreme value theory once the
tails of the $X_n$ is known, which is precisely what has been
the objective of this paper.

Whether this asymptotic scheme is the most relevant one is, however,
questionable. One could equally well assume the job length fixed
at $T$ and the length of the $N$ subjobs to be $T/N$, and intermediate
possibilities. This leads into specific questions on extreme value
theory in a triangular array setting, which are currently
under investigation.

\medskip

\noindent {\bf C)} An alternative to the Cram\'er-Lundberg
theory for geometric sums that has been one of our main tools
is what could suitably be called Renyi theory, cf.\ \cite{Volodya}.
One considers there
a weak convergence triangular setting where still $x\to\infty$
but the parameters of the geometric sum depend on $x$;
this is also related to the heavy-traffic or diffusion limit
setting of risk and queueing theory, cf.\ e.g.\ \cite{RP} and
\cite{APQ} X.7. 
Renyi theory (e.g.\ \cite{Volodya}) provides the following alternative to 
Proposition \ref{Prop20.10a}:
\begin{Prop}\label{Prop20.10axq} 
For each $x$, let $V_1(z),V_2(z),\ldots$ be i.i.d.\ 
with common density $k(v;z)$ and
$N(z)\in\nat$ an independent r.v.\ with $\Prob(N(z)=n)=(1-\rho(z))\rho(z)^n$,
$S(z)=V_1(z)+\cdots+V_N(z)$. If $\rho(z)\to 1$ as $z\to\infty$ 
and $V_k(z)\convdistr V$, $\Exp V_k(z)\to 1/\mu$ for some r.v.\
$V$ with finite mean $\mu^{-1}$, then $\mu(1-\rho(z)S(z)$ has a limiting
standard exponential distribution.
\end{Prop}

\begin{Cor}\label{Cor03.01a} In the RESTART setting,
$\mu\overline G(t)S(t)$ has a limiting
standard exponential distribution as $t\to\infty$.
\end{Cor} 
The implication is that
\begin{equation}\label{03.01a}
\Prob\bigl(\mu\overline G(t)S(t)>y\bigr)\,\to\,\e^{-y}
\end{equation}
for any fixed $y$. Noting that $\mu\overline G(t)\sim\gamma(t)$
and replacing $y$ by $\gamma(t)x$, this suggests $\Prob(S(t)>x)$
$\approx$ $\e^{-\gamma(t)x}$, i.e.\ the Cram\'er-Lundberg approximation.
Of course, the derivation is not rigorous since \refs{03.01a}
requires that $y$ is fixed. Nevertheless, it is indeed possible to derive
some of our results from \refs{03.01a}. The main reason that we have chosen
Cram\'er-Lundberg asymptotics as our basic vehicle is that simple
bounds are available (Lundberg's inequality) which is not the case
for Renyi theory.

\medskip

\noindent {\bf D)} Since $H$ is a mixture of the $H_t$ given by 
$H_t(x)=\Prob\bigl(S(t)+t\le x\bigr)$
and the tail  of $H_t$ obeys the Cram\'er-Lundberg asymptotics,
we are dealing with the problem of determining the tail of a mixture
where the tails of the mixing components are known. Looking for
literature on this problem, we found a set of papers emerging from
reliability and survival analysis 
( Finkelstein \& Esaulova~\cite{Fink} and references there)
which suggest our logarithmic
asymptotics results but do not prove them because the assumptions
are too stringent to apply to our setting.

\medskip

\noindent {\bf E)} The $\approx_{\log\log}$
asymptotics in part (1:2) identifies $\theta_{12}$ as the correct
exponent to $x$ in $\log \overline H(x)$, but does not allow
sharpenings like $  \overline H(x)\logsim\e^{-c_{12}x^{\theta_{12}}}$,
$  \overline H(x)\logsim\e^{-c_{12}\log^{q_{12}}x\, x^{\theta_{12}}}$
etc. Inspection of the proof shows that to 
obtain such strengthenings, one needs first of all to 
be able to replace
the $b$ in Lemma~\ref{L21.10aaa} with a fixed value $\beta$ rather
than considering $b$'s arbitrarily close to $\beta$. This would be
the case if, e.g., one assumed $\overline G$ to be regularly varying
with index $-\beta$ rather than just $\overline G(t)\logsim t^{-\beta}$.
This does not appear to be all that restrictive, but does not suffice
since one also needs to replace the $\pm\epsilon$ in  Lemma~\ref{L20.10c}
with sharper bounds. This amounts to second-order asymptotics of
the $\overline H_t(x)$, i.e.\ to obtain  second-order uniform
Cram\'er-Lundberg expansions which does not appear easy at all.

\end{document}